\documentclass[12pt]{article}

\usepackage{amssymb}
\usepackage{amsmath}

\def\disp{\displaystyle}
\def\o{\over}
\def\T{{\cal T}}
\def\reals{I\!\!R}
\def\a{\alpha}
\def\ve{\varepsilon}
\def\g{\gamma}
\def\gm{\gamma(dy,du)}
\def\G{\Gamma}
\def\U{{\cal U}}
\def\P{{\cal P}}

\def\ph{\varphi}
\def\l{\lambda}
\def\d{\delta}
\def\t{\tau}
\def\bu{\bar u}
\def\by{\bar y}
\def\ve{\varepsilon}
\def\k{\kappa}
\def\b{\beta}

\def\up{\uparrow}

\def\s{\sigma}
\def\ts{\tilde\sigma}
\def\hf{\hfill{$\Box$}}
\def\D{\Delta}

\newtheorem{Theorem}{Theorem}[section]
\newtheorem{Proposition}[Theorem]{Proposition}
\newtheorem{Remark}[Theorem]{Remark}
\newtheorem{Lemma}[Theorem]{Lemma}
\newtheorem{Corollary}[Theorem]{Corollary}

\begin{document}

{\bf \Large{Linear Programming Formulations of Deterministic Infinite Horizon Optimal Control Problems in Discrete Time}}
 \bigskip
 \bigskip
 \bigskip

 {\bf V. Gaitsgory$^a$, A. Parkinson$^a$ and I. Shvartsman$^b$}\\
 $^a$ {\it\small{Department of Mathematics, Macquarie University, Eastern Road, Macquarie Park, NSW 2113, Australia } }\\
 $^b$ {\it \small{Department of Mathematics and Computer Science,
Penn State Harrisburg, Middletown, PA 17057, USA}}

 \bigskip
 \bigskip
 \bigskip

{\bf Abstract.} \small{This paper is devoted to a study of  infinite horizon optimal control problems with time discounting and time averaging criteria in discrete time.
We establish that these problems are related to certain  infinite-dimensional linear programming  (IDLP) problems. We also establish asymptotic relationships between the optimal values of problems with time discounting and long-run average criteria.}

\bigskip


\bigskip
 \bigskip
 \bigskip

\section{Introduction}

The linear programming (LP) approach to control systems is based on the fact that the
occupational measures generated by admissible controls and the corresponding solutions of
a dynamical system satisfy certain linear equations that represent the system's dynamics in an
integral form. The idea of such linearization was explored extensively in both deterministic and stochastic settings (see, e.g.,  \cite{BhBo}, \cite{Vivek}, \cite{BGQ}, \cite{F-V}, \cite{Kurtz},  \cite{Stockbridge}, \cite{Stockbridge1}   and, respectively, \cite{Adelman}, \cite{Gai-F-Leb}, \cite{Gai8}, \cite{GQ}, \cite{GQ-1}, \cite{GR}, \cite{Goreac-Serea}, \cite{Her-Her-Lasserre}, \cite{Adelman-1}, \cite{Lass-Trelat}, \cite{QS}, \cite{Rubio}, \cite{Vinter}     as well as references therein). In \cite{GQ} and \cite{GQ-1}  in particular, the validity of  LP formulations of deterministic infinite time horizon problems of optimal control
with time average and  time discounting criteria was proved for
systems evolving in continuous time (note that other approachers/techniques for dealing with deterministic optimal control problems on the infinite time horizon have been studied, e.g., in   \cite{BCD}, \cite{B12}, \cite{CHL}, \cite{Z14}; see also references therein).
 In the present paper, we show that the LP formulations of
 problems of optimal control
with time average and  time discounting criteria are valid for systems evolving in discrete time.

Note that some of the results of \cite{GQ} and \cite{GQ-1} were obtained under certain technical assumptions. For example, the statement implying the validity of the LP formulation of the long run average optimal control problem (see Theorem 2.6 in \cite{GQ-1}) was proved under the assumption that the dependence of the control set on the state variables is Lipschitz continuous. These assumptions can be significantly relaxed
in dealing with the discrete time systems. In particular, the result about the validity of the LP formulation
of the long run average optimal control problem in  discrete time   is established in this paper
under the assumption that the dependence of the control set on the state variables is upper semicontinuous.
Also, it is worth noting that the results in  \cite{GQ-1} (see also Remark 4.5 in \cite{GQ}) are stated with the use of  the relaxed controls formalism, the latter playing no role in tackling the discrete time systems.

Everywhere in what follows we will be dealing with the discrete time controlled dynamical system
\begin{equation}\label{A1}
\begin{aligned}
&y(t+1)=f(y(t),u(t)), \; t=0,1,\dots\, \\
&y(0)=y_0,\\
&y(t)\in Y,\\
&u(t)\in U(y(t)).
\end{aligned}
\end{equation}
Here  $Y$ is a given nonempty compact subset of $\reals^m$, $\ U(\cdot):\,Y\leadsto U_0$ is an upper semicontinuous compact-valued mapping to a given compact metric space $U_0$,
$\ f(\cdot, \cdot ):\,\reals^m\times U_0\to \reals^m$ is a continuous function.

Note that the last two constraints of \eqref{A1} can be rewritten as one:
$$
u(t)\in A(y(t)),
$$
where the map $\ A(\cdot):\,Y\leadsto U_0$ is defined by the equation
\begin{equation*}
\begin{aligned}
A(y):=\{u\in U(y)|\, f(y,u)\in Y\} \ \ \ \forall y\in Y.
\end{aligned}
\end{equation*}
As can be readily verified, the map $A(\cdot)$ is upper semicontinuous and its graph $G$,
$$
G:={\rm graph}\,A=\{(y,u)|\,y\in Y,\,u\in U(y),\,f(y,u)\in Y\},
$$
is a compact subset of $Y\times U_0$.

A control $u(\cdot)$ and the pair $(y(\cdot),u(\cdot))$ will be called an admissible control and, respectively, an admissible process if the relationships \eqref{A1} are satisfied.
The sets of admissible controls will be denoted by $\U(y_0)$ or $\U_S(y_0)$, depending on whether the problem is considered on the  infinite  time horizon ($t\in \T:=\{0,1,\dots\}$) or on a finite time sequence  ($t\in \{0,\dots,S-1\}$, where $S$ is a positive integer).

Consider the  optimal control problem
\begin{equation}\label{A111}
\min_{u(\cdot)\in \U(y_0)}\sum_{t=0}^{\infty} \a^t g(y(t),u(t))=:V_{\a}(y_0),
\end{equation}
where $g:\,\reals^m\times U_0\to \reals^m$ is a continuous function and  $\a\in (0,1)$ is a discount factor.
Consider also the optimal control problem
\begin{equation}\label{A112-1}
\min_{u(\cdot)\in \U_S(y_0)}\sum_{t=0}^{S-1} g(y(t),u(t))=:V(S, y_0),
\end{equation}
 Everywhere in the paper, it is assumed that

A1. {\em The set $\U(y_0)$ is not empty (that is, there exists at least one admissible control).}

As shown below (see Propositions \ref{P1} and \ref{Prop-exist-ave}), the minima in (\ref{A111}) and
(\ref{A112-1}) are achieved if A1 is satisfied. To obtain our main results, we  use a stronger assumption:

A2. {\em The set $A(y)$ is not empty for any $y\in Y$.}

This assumption implies non-emptiness of $\U(y)$ for any $y\in Y$  (systems that satisfy such a property are called {\em viable}; see \cite{Aub}).

Along with optimal control problems (\ref{A111}) and  (\ref{A112-1}), let us consider  two infinite-dimensional (ID) linear programming (LP) problems:
\begin{equation}\label{D1}
\min_{\g\in W_{\a}(y_0)} \int_{G} g(y,u)\g(dy,du):=g^*_{\a}(y_0)
\end{equation}
and
\begin{equation}\label{M22}
\min_{\g\in W} \int_{G} g(y,u)\g(dy,du):=g^*,
\end{equation}
where
$W_{\a}(y_0)$ and $W$ are  subsets of $\P(G)$  (here and in what follows $\P(G)$ stands for the space of probability measures on Borel subsets of $G$) defined by the equations:
 \begin{equation}\label{e-W}
\begin{aligned}
W_{\a}(y_0):=\{&\g\in \P(G)|\, \\
&\int_{G}[\a(\ph(f(y,u))-\ph(y))
+(1-\a)(\ph(y_0)-\ph(y))]\gm=0\quad \forall \ph\in C(Y)\}
\end{aligned}
\end{equation}
and
\begin{equation}\label{M17}
\begin{aligned}
W:=\{\g\in \P(G)|\, \int_{G}(\ph(f(y,u))-\ph(y))
\gm=0\quad \forall \ph\in C(Y)\}.
\end{aligned}
\end{equation}
Note that (\ref{D1}) and (\ref{M22}) are indeed LP problems since both the objective functions and the constraints defining  $W_{\a}(y_0)$ and $W$ are linear in the \lq\lq decision variable" $\g$. Note also  that $W$ can be obtained from $W_{\a}(y_0)$ by setting $\a=1$.

In the paper, we prove that,
under Assumption A2,
\begin{equation}\label{eq-Res1}
(1-\a)V_{\a}(y_0)=g^*_{\a}(y_0)
\end{equation}
 and  the limits $\ \lim_{\a\up 1}\min_{y\in Y}(1-\a)V_{\a}(y) $
 and $\ \lim_{S\to \infty}\min_{y\in Y}{1\o S}V(S,y) $ exist and are equal to $g^*$:
\begin{equation}\label{eq-Res2}
\lim_{\a\up 1}\min_{y\in Y}(1-\a)V_{\a}(y)=\lim_{S\to \infty}\min_{y\in Y}{1\o S}V(S,y)=g^*.
\end{equation}
It is worth mentioning that there exists an extensive literature devoted to the relationship between the limits of the sums $\ \disp (1-\a)\sum_{t=0}^{\infty} \a^t b_t$ and
$\ \disp {1\o S}\sum_{t=0}^{S-1}b_t$ as $\a\up 1$ and $S\to \infty$, respectively. There are many examples showing that these limit may not exist (see, e.g., \cite{Bishop-Feinberg}, where relationships between the corresponding lower and upper limits were investigated). However, provided that the sequence $\{b_t\}$ is bounded, the existence of one of these limits implies the existence of the other and their equality (see, e.g., \cite{Filar92}).
In the context of optimal control in discrete time, relationships between the lower and upper limits of
$(1-\a)V_{\a}(y) $ and ${1\o S}V(S,y) $
 were studied, e.g., in \cite{Sorin92} and \cite{Renault09}. The (full) aforementioned limits may not exist, and, as was shown in \cite{Sorin92} (without the assumption about the compactness of the set of admissible states $Y$),  these limits, even if exist, may be different.
As mentioned above, in this paper we establish that, under the validity of A2, the limits of the minima over the initial conditions of $(1-\a)V_{\a}(y) $ and ${1\o S}V(S,y) $
exist and are equal to the optimal value of the IDLP problem (\ref{M22}).


 The paper is organized as follows. Section \ref{Preliminaries} contains some preliminary results used in the sequel. In Section \ref{Sec-occup-meas}, we introduce  discounted and \lq\lq non-discounted" occupational measures and we reformulate problems (\ref{A111}) and (\ref{A112-1}) in terms of minimization over the sets of such measures. In Section \ref{Sec-equality}, we establish that (\ref{eq-Res1}) is valid, and in Section \ref{Sec-asymptotic} we prove
the validity of (\ref{eq-Res2}). In this section,  we also establish asymptotic properties of the sets of discounted and non-discounted occupational measures. In Section \ref{Sec-Appendix}, we prove auxiliary results that are used in Sections \ref{Sec-equality} and \ref{Sec-asymptotic}.

\section{Preliminaries}\label{Preliminaries}

Everywhere in this and the following sections, it is assumed that A1 is satisfied.

\begin{Proposition}\label{P1}
The minimum in \eqref{A111} is achieved.
\end{Proposition}

{\bf Proof.} For an admissible process $(y(\cdot),u(\cdot))$, denote $J_{\a}(u,y_0):=\sum_{t=0}^{\infty} \a^t g(y(t),u(t))$. Let $u_k(\cdot)$, $k=1,2,\dots$ be a minimizing sequence of controls and let $y_k(\cdot)$ be the corresponding sequence of trajectories. By using the diagonalization argument and taking into account compactness of $G$, we can find convergent subsequences (we do not relabel) $u_k(t)\to \bu(t)$ and $y_k(t)\to \by(t)$ for all $t$. By passing to the limit in the relation $y_k(t+1)=f(y_k(t),u_k(t))$ as $k\to \infty$ we conclude that the process $(\by(\cdot),\bu(\cdot))$ is admissible.
For any natural $N$ we have
$$
|J_{\a}(u_k,y_0)-J_{\a}(\bu,y_0)|\le \sum_{t=0}^{N}\a^t |g(y_k(t),u_k(t))-g(\by(t),\bu(t))|+\sum_{t=N+1}^{\infty}\a^t |g(y_k(t),u_k(t))-g(\by(t),\bu(t))|.
$$
Take $\ve>0$ and find $N$ large enough so that the second sum does not exceed $\ve/2$ for all $k$, then the first sum can be made less than $\ve/2$ by taking sufficiently large $k$.
Therefore, $J_{\a}(u_k,y_0)\to J_{\a}(\bu,y_0)$ as $k\to \infty$, which implies that the process $(\by(\cdot),\bu(\cdot))$ is optimal.


\begin{Proposition}\label{P2}
The optimal value function $V_{\a}(\cdot)$ is lower semicontinuous.
\end{Proposition}

{\bf Proof.} Take a sequence $y_{0k}\to y_0$ as $k\to \infty$ such that $V_{\a}(y_{0k})<\infty$. Let $u_k(\cdot)$ be the corresponding sequence of minimizing controls, that is, controls such that $V_{\a}(y_{0k})=J_{\a}(u_k,y_{0k})$. We want to show that $\disp \liminf_{k\to \infty} V_{\a}(y_{0k})\ge V_{\a}(y_0).$
Without loss of generality assume that $\disp \liminf_{k\to \infty} V_{\a}(y_{0k})$ is reached on the same sequence $y_{0k}$. Again, using the diagonalization argument and passing to a subsequence, we can assume that $u_k(t)$ converges to admissible control $\bu(t)$ for all $t$. Using the same argument as in the proof of Proposition \ref{P1} we can show that
$\lim_{k\to \infty} J_{\a}(u_k,y_{0k})= J_{\a}(\bu,y_0)$. We have
$$
\lim_{k\to \infty} V_{\a}(y_{0k})=\lim_{k\to \infty} J_{\a}(u_k,y_{0k})= J_{\a}(\bu,y_0)\ge \min_{u(\cdot)\in \U(y_0)} J_{\a}(u,y_0)=V_{\a}(y_0),
$$
which is the required inequality. \hfill{$\Box$}

\begin{Proposition}\label{Prop-exist-ave}
The minimum in (\ref{A112-1}) is achieved and the optimal value function $V(S,\cdot)$ is lower semicontinuous.
\end{Proposition}

{\bf Proof.} The fact that the minimum in (\ref{A112-1}) is  achieved is obvious  (since it is a finite-dimensional problem on a compact set), and the fact that $V(S,\cdot)$ is lower semicontinuous is proved similarly to Proposition \ref{P2}.
\hfill{$\Box$}

\begin{Corollary}\label{Cor-min-achieved}
The minima in (\ref{eq-Res2}) are achieved.
\end{Corollary}

{\bf Proof.} The proof follows from the fact that the functions $V_{\alpha}(\cdot) $ and $V(S,\cdot) $ are lower semicontinuous.
\hfill{$\Box$}

\begin{Proposition}\label{P3}
For any $y\in Y$ such that $V_{\a}(y)<\infty$, the following equation is valid
\begin{equation}\label{B1}
V_{\a}(y)=\min_{u\in A(y)}\{g(y,u)+\a V_{\a}(f(y,u))\}.
\end{equation}
\end{Proposition}

{\bf Proof.} The proposition is the well known dynamic programming principle for problem \eqref{A111}. For completeness of the exposition, we reproduce its proof in Section \ref{Sec-Appendix}. \hfill{$\Box$}

For a lower semicontinuous function $\psi:\,Y\to \reals$, let $H_{\psi}(y)$ be defined as follows
\begin{equation*}
\begin{aligned}
H_{\psi}(y):=\min_{u\in A(y)}\{\a (\psi(f(y,u))-\psi(y))+g(y,u)\}.
\end{aligned}
\end{equation*}
Then equation \eqref{B1} can be written as
\begin{equation}\label{E3}
H_{V_{\a}}(y)-(1-\a)V_{\a}(y)=0,
\end{equation}
which resembles the Hamilton-Jacobi-Bellman equation for continuous time systems; see, e.g., \cite{BCD}.

\section{Occupational Measure Formulations}\label{Sec-occup-meas}

Let $(y(\cdot), u(\cdot))$ be an admissible process. A probability measure $\g^{\a}_{(y(\cdot),u(\cdot))}$ is called the {
\em discounted occupational measure} generated by the process $(y(\cdot), u(\cdot))$ if, for any Borel set $Q\subset G$,
\begin{equation}\label{E6}
\g^{\a}_{(y(\cdot),u(\cdot))}(Q)=(1-\a) \sum_{t=0}^{\infty} \a^t 1_Q(y(t),u(t)),
\end{equation}
where $1_Q(\cdot)$ is the indicator function of $Q$. A probability measure
$\g_{(y(\cdot), u(\cdot)),S}$ is called the {\em  occupational measure} generated by the process $(y(\cdot), u(\cdot))$ over the time sequence $\{0,1,...,S-1 \}$ if, for any Borel set $Q\subset G$,
\begin{equation*}
\g_{(y(\cdot), u(\cdot)),S}(Q)={1\o S} \sum_{t=0}^{S-1} 1_Q(y(t),u(t)),
\end{equation*}

It can be shown that if $\g^{\a}_{(y(\cdot),u(\cdot))}$ is the discounted occupational measure generated by the process $(y(\cdot), u(\cdot))$, then
\begin{equation}\label{G8}
\int_{G} q(y,u) \g^{\a}_{(y(\cdot),u(\cdot))}(dy,du)=(1-\a)\sum_{t=0}^{\infty} \a^t q(y(t),u(t))
\end{equation}
for any Borel measurable function $q$ on $G$. Also, it can be shown that if $\g_{(y(\cdot), u(\cdot)),S}$ is the occupational measure generated by the process $(y(\cdot), u(\cdot))$ over the time sequence $\{0,1,...,S-1 \}$, then
\begin{equation}\label{G88}
\int_{G} q(y,u) \g_{(y(\cdot),u(\cdot)),S}(dy,du)={1\o S}\sum_{t=0}^{S-1} q(y(t),u(t))
\end{equation}
for any Borel measurable function $q$ on $G$.

To describe convergence properties of occupational measures, we introduce the following metric on $\P(G)$:
$$
\rho(\g',\g''):=\sum_{j=1}^{\infty} {1\o 2^j}\left|\int_G q_j(y,u)\g'(dy,du)-\int_G q_j(y,u)\g''(dy,du)\right|
$$
for $\g',\g''\in \P(G)$, where $q_j(\cdot),\,j=1,2,\dots,$ is a sequence of Lipschitz continuous functions dense in the unit ball of the space of continuous functions $C(G)$ from $G$ to $\reals$.
This metric is consistent with the weak$^*$ convergence topology on $\P(G)$, that is,
a sequence $\g^k\in \P(G)$ converges to $\g\in \P(G)$ in this metric if and only if
$$
\lim_{k\to \infty}\int_G q(y,u)\g^k(dy,du)=\int_G q(y,u)\g(dy,du)
$$
for any $q\in C(G)$. Note that the sets $W_{\alpha}(y_0)$ and $W$ are compact in this topology.

Using the metric $\rho$, we can define the ``distance" $\rho(\g,\Gamma)$ between $\g\in \P(G)$ and $\Gamma\subset \P(G)$
and the Hausdorff metric $\rho_H(\Gamma_1,\Gamma_2)$ between $\Gamma_1\subset \P(G)$ and $\Gamma_2\subset \P(G)$ as follows:
$$
\rho(\g,\Gamma):=\inf_{\g'\in \Gamma}\rho(\g,\g'),\quad
\rho_H(\Gamma_1,\Gamma_2):=\max\{\sup_{\g\in \Gamma_1}\rho(\g,\Gamma),\sup_{\g\in \Gamma_2}\rho(\g,\Gamma_2)\}.
$$
Note that, although, by some abuse of terminology,  we refer to
$\rho_H(\cdot,\cdot)$ as  a metric on the set of subsets of
${\cal P} (Y \times U)$, it is, in fact, a semi metric on this set
(since $\rho_H(\Gamma_1, \Gamma_2)=0$ implies  $\Gamma_1
= \Gamma_2$ if  $\Gamma_1$ and $\Gamma_2$ are closed and the equality may not be true if at least one of these sets is not closed).

Introduce the following notation for the sets of occupational measures:
$$
\G_{\a}(y_0):=\bigcup_{_{u(\cdot)\in \U(y_0)}}\{\g^{\a}_{(y(\cdot),u(\cdot))}\},\quad \ \ \ \ \G_{\a}:=\bigcup_{y_0\in Y}\{\G_{\a}(y_0)\},
$$
$$
\G(S,y_0):=\bigcup_{u(\cdot)\in \U_S(y_0)}\{\g_{(y(\cdot),u(\cdot)),S}\},\quad \ \ \ \ \G(S):=\bigcup_{y_0\in Y}\{\G_S(y_0)\}.
$$
Due to \eqref{G8} and \eqref{G88}, problems \eqref{A111} and (\ref{A112-1}) can be rewritten in the form
\begin{equation}\label{A3}
\min_{\g\in \G_{\a}(y_0)} \int_{G} g(y,u)\g(dy,du)=(1-\a)V_{\a}(y_0)
\end{equation}
and
\begin{equation*}
\min_{\g\in \G(S,y_0)} \int_{G} g(y,u)\g(dy,du)=V(S,y_0),
\end{equation*}
respectively.

\section{Validity of (\ref{eq-Res1})}\label{Sec-equality}

\begin{Proposition}\label{P4}
The inclusion $\G_{\a}(y_0)\subset W_{\a}(y_0)$ is true.
\end{Proposition}
{\bf Proof.} For arbitrary $\ph\in C(Y)$ and admissible process $(y(\cdot),u(\cdot))$  we have
\begin{equation*}
\sum_{t=0}^{\infty} \a^t \ph(y(t))=\ph(y_0)+\a\sum_{t=0}^{\infty} \a^t \ph(f(y(t),u(t))).
\end{equation*}
Multiplying both sides by $1-\a$ and taking into account \eqref{G8}, we obtain
$$
\int_{G}\ph(y)\,\gamma^{\a}_{(y(\cdot),u(\cdot))}(dy,du)=\int_{G}(1-\a)\ph(y_0)\,\gamma^{\a}_{(y(\cdot),u(\cdot))}(dy,du)+\a \int_{G}\ph(f(y,u))\,\gamma^{\a}_{(y(\cdot),u(\cdot))}(dy,du),
$$
where $\g^{\a}_{(y(\cdot),u(\cdot))}\in \G_{\a}(y_0)$ is generated by $(y(\cdot),u(\cdot))$. The latter is equivalent to
\begin{equation*}\label{C1}
\int_{G} [\a(\ph(f(y,u))-\ph(y))+(1-\a)(\ph(y_0)-\ph(y))]\gamma^{\a}_{(y(\cdot),u(\cdot))}(dy,du)=0.
\end{equation*}
This implies that $\g^{\a}_u\in W_{\a}(y_0)$, which concludes the proof of the proposition.
\hfill{$\Box$}

\begin{Remark}
Due to the assumed validity of A1,  $\G_{\a}(y_0)\neq \emptyset$ and, hence, $W_{\a}(y_0)\neq \emptyset$.
\end{Remark}
Note that from Proposition \ref{P4} it follows that
\begin{equation}\label{A5}
g^*_{\a}(y_0)\le (1-\a)V_{\a}(y_0).
\end{equation}

\bigskip

Let $LS$ be the class of bounded lower semicontinuous functions from $Y$ to $\reals$.
Note that $V_{\a}(\cdot)\in LS$ if Assumption A2 is satisfied. In fact, in this case
\begin{equation}\label{e-bounded}
\left|V_{\a}(y)\right|\ \leq \ \frac{M}{1-\a}\ \ \forall y\in Y, \ \ \ \ \ {\rm where} \ \ \ \ \  M:=\max_{(y,u)\in G}|g(y,u)|.
\end{equation}
From this point on, it is everywhere assumed that Assumption A2 is indeed satisfied.

Consider the max-min problem
\begin{equation}\label{B5}
\begin{aligned}
&\sup_{\psi\in LS}\inf_{y\in Y}\{H_{\psi}(y)+(1-\a)(\psi(y_0)-\psi(y))\}=\\
&\sup_{\psi\in LS}\inf_{(y,u)\in G}\{g(y,u)+\a(\psi(f(y,u))-\psi(y))+(1-\a)(\psi(y_0)-\psi(y))\}:=\mu^*_{\a}(y_0).
\end{aligned}
\end{equation}

We say that $\tilde \psi$ is a solution of \eqref{B5} if
$$
\inf_{(y,u)\in G}\{g(y,u)+\a(\tilde\psi(f(y,u))-\tilde\psi(y))+(1-\a)(\tilde\psi(y_0)-\tilde\psi(y))\}=\mu^*(y_0).
$$
Our first main result is the following theorem.
\begin{Theorem}\label{P5}
The optimal values in problems \eqref{D1} and \eqref{B5} coincide and are equal to the optimal value of \eqref{A111} multiplied by $(1-\a)$, that is,
\begin{equation}\label{D3}
\mu^*_{\a}(y_0)=g^*_{\a}(y_0)=(1-\a)V_{\a}(y_0).
\end{equation}
Moreover, the supremum in \eqref{B5} is reached at $\psi=V_{\a}$.
\end{Theorem}

{\bf Proof.} From Proposition \ref{P3} we have
\begin{equation*}
\min_{u\in A(y)}\{g(y,u)+\a V_{\a}(f(y,u))-V_{\a}(y)\}=0 \quad\hbox{for all }y\in Y,
\end{equation*}
which implies that
\begin{equation*}
\min_{(y,u)\in G}\{g(y,u)+\a V_{\a}(f(y,u))-V_{\a}(y)\}=0.
\end{equation*}
Therefore,
\begin{equation}\label{D2}
\begin{aligned}
(1-\a)V_{\a}(y_0)= \min_{(y,u)\in G} \{g(y,u)+\a(V_{\a}(f(y,u))-V_{\a}(y))+(1-\a)(V_{\a}(y_0)-V_{\a}(y))\}
\le \mu^*_{\a}(y_0).
\end{aligned}
\end{equation}
Taking into account \eqref{A5}, we get
\begin{equation}\label{A20}
g^*_{\a}(y_0)\le (1-\a)V_{\a}(y_0)\le \mu^*_{\a}(y_0).
\end{equation}
Let us show the opposite inequality. For $\psi\in LS$ denote
\begin{equation}\label{A7}
\mu_{\a}(\psi,y_0):=\inf_{(y,u)\in G}\{g(y,u)+\a (\psi(f(y,u))-\psi(y))+(1-\a)(\psi(y_0)-\psi(y))\},
\end{equation}
so that $\disp \mu^*_{\a}(y_0)=\sup_{\psi\in LS}\mu_{\a}(\psi,y_0)$. Take $\g\in W_{\a}(y_0)$, arbitrary $\psi\in LS$ and let $\{\psi_n\}_{n=1}^{\infty}$ be a bounded sequence of continuous  functions such that $\psi_n(y)\to \psi(y)$ point-wise on $Y$ as $n\to \infty$ (due to (\ref{e-bounded}), such a sequence exists; see, e.g., Theorem A6.6 in \cite{Ash}). From \eqref{A7}, from  Lebesgue dominated convergence theorem and from the definition of $W_{\a}(y_0)$ it follows that
\begin{equation*}
\begin{aligned}
\mu_{\a}(\psi,y_0)\le& \int_{G}[g(y,u)+\a (\psi(f(y,u))-\psi(y))+(1-\a)(\psi(y_0)-\psi(y))]\gm\\
=&\lim_{n\to \infty} \int_{G}[g(y,u)+\a (\psi_n(f(y,u))-\psi_n(y))+(1-\a)(\psi_n(y_0)-\psi_n(y))]\gm\\
=&\int_{G}g(y,u)\gm.
\end{aligned}
\end{equation*}
Taking supremum with respect to $\psi\in LS$ and minimum with respect $\g\in W_{\a}(y_0)$ leads to $\mu^*_{\a}(y_0)\le g^*_{\a}(y_0)$ which, together with
\eqref{A20}, implies \eqref{D3}. It also follows from \eqref{D2} that
$$
\mu^*_{\a}(y_0)=\min_{(y,u)\in G} \{g(y,u)+\a(V_{\a}(f(y,u))-V_{\a}(y))+(1-\a)(V_{\a}(y_0)-V_{\a}(y))\},
$$
which implies the second part of the theorem.
\hfill{$\Box$}

\begin{Corollary}
The following equality is valid
\begin{equation}\label{e-disc-equality}
\bar{\rm co}\,\G_{\a}(y_0)=W_{\a}(y_0),
\end{equation}
where $\bar{\rm co}\, $ stands for the closure of the convex hull of the corresponding set.
\end{Corollary}
{\bf Proof.} Due to (\ref{D1}) and (\ref{A3}), the equality (\ref{eq-Res1}) can be rewritten in the form
$$
\min_{\g\in \G_{\a}(y_0)} \int_{G} g(y,u)\g(dy,du) = \min_{\g\in W_{\a}(y_0)} \int_{G} g(y,u)\g(dy,du),
$$
which implies that
$$
\min_{\g\in \bar{\rm co}\,\G_{\a}(y_0)} \int_{G} g(y,u)\g(dy,du) = \min_{\g\in W_{\a}(y_0)} \int_{G} g(y,u)\g(dy,du).
$$
Since the latter is valid for any continuous $g$, it proves the validity of (\ref{e-disc-equality}).
\hfill{$\Box$}
\begin{Remark}
{\rm Note that problem (\ref{B5}) can be shown to be equivalent to the problem dual to the IDLP problem (\ref{D1}) (see Appendix of \cite{GQ}), with the equality of the optimal values being a part of the duality relationships between these two problems.}
\end{Remark}

\section{Validity of (\ref{eq-Res2})}\label{Sec-asymptotic}

 Let us introduce the following notation:

\begin{equation}\label{A112}
{1\o S}\min_{y_0\in Y}\min_{u(\cdot)\in \U_S(y_0)}\sum_{t=0}^{S-1} g(y(t),u(t))=:G_S,
\end{equation}
where the minimization is over admissible controls  and over the initial conditions in $Y$.

The main results of this section are Theorems \ref{T1} and \ref{T2} below. In Theorem \ref{T1} we, in particular, establish existence and equality of the limits in \eqref{eq-Res2}.
Theorem \ref{T2} deals with a limiting property of the sets of occupational measures and is closely related to Theorem \ref{T1}.
Continuous-time analogs of Theorems \ref{T1} and \ref{T2} are proved in \cite{GQ}, Chapter 6. However, in continuous time, as opposed to discrete time, a few strong assumptions are needed for the validity of the corresponding results (e.g., Lipschitz continuity of the value function).

Let

\begin{equation}\label{M21}
\mu^*:=\sup_{\psi\in LS}\inf_{(y,u)\in G}\{g(y,u)+\psi(f(y,u))-\psi(y)\}.
\end{equation}

\begin{Theorem}\label{T1}
The limits $\disp \lim_{\a\up 1}\min_{y\in Y}(1-\a)V_{\a}(y)$ and $\disp\lim_{S\to \infty}G_S$ exist and
$$
\lim_{\a\up 1}\min_{y\in Y}(1-\a)V_{\a}(y)=\lim_{S\to \infty}G_S=g^*=\mu^*.
$$
\end{Theorem}

The proof is broken down into a series of propositions and lemmas.

\begin{Proposition}\label{P10}
The equality $g^*=\mu^*$ holds true.
\end{Proposition}
{\bf Proof.}
Take any $\psi\in LS$. Integrating the inequality
$$
g(y,u)+\psi(f(y,u))-\psi(y) \ge \inf_{(y,u)\in G}\{g(y,u)+\psi(f(y,u))-\psi(y)\}
$$
with respect to arbitrary $\g\in W$ we obtain
$$
\int_{G} g(y,u)\g(dy,du)\ge \inf_{(y,u)\in G}\{g(y,u)+\psi(f(y,u))-\psi(y)\}.
$$
Taking minimum with respect to $\g\in W$ and supremum with respect to $\psi\in LS$, we conclude that
\begin{equation}\label{M26}
g^*\ge \mu^*.
\end{equation}
Let us show the opposite inequality.
Define
\begin{equation}\label{D10}
\mu_C^*:=\sup_{\psi\in C(Y)}\min_{(y,u)\in G}\{g(y,u)+\psi(f(y,u))-\psi(y)\},
\end{equation}
that is, compared to \eqref{M21}, supremum in the formula above is taken with respect to continuous, rather than lower semicontinuous bounded functions.
It is clear that
\begin{equation}\label{M24}
\mu_C^*\le \mu^*,
\end{equation}
therefore $\mu_C^*<\infty$.

Let $\{\phi_i\}_{i=1}^{\infty}$ be a sequence of functions in $C(Y)$ with the following properties: (i) any finite collection of functions from this sequence is linearly independent on $Y$, (ii) for any $\psi\in C(Y)$ and any $\d>0$ there exist $N$ and scalars $\l_i^N$, $i=1,\dots,N$ such that $\disp \sup_{y\in Y}|\psi(y)-\sum_{i=1}^N \l_i^N\phi_i(y)|\le \d$. (An example of such sequence is the sequence of monomials $y_1^{i_1}\dots y_m^{i_m},\,i_1,\dots,i_m=0,1,\dots$, where $y_j$ stands for the $j$th component of $y$.)

Let us notice first that for any $\psi\in C(Y)$ we have
\begin{equation}\label{H3}
\min_{(y,u)\in G}\{\psi(f(y,u))-\psi(y)\}\le 0.
\end{equation}
Indeed, if this was not the case, then, for $\psi_m:=m\psi$ with positive integer $m$ we would get
$$
\lim_{m\to \infty}\min_{(y,u)\in G}\{g(y,u)+\psi_m(f(y,u))-\psi_m(y)\}=+\infty,
$$
which contradicts boundedness of $\mu_C^*$.

Assume that functions $\{\phi_i\}$ are normalized so that $\max_{y\in Y}|\phi_i(y)|<1/2^i$. Define $\hat{Q}\subset\reals\times l^1$ by
		\begin{equation*}
				\begin{aligned}
			\hat{Q} := &\{(\theta,x) |\, \theta\geq\int_{G}g(y,u)\gamma(dy,du),\, x=\left(x_1,x_2,\dots\right),  \\
			&x_i=\int_{G}(\phi_i\left(f\left(y,u\right)\right)-\phi_i\left(y\right)))\gamma\left(dy,du\right),\,\gamma\in\mathcal{P}(G)\}.
		\end{aligned}
		\end{equation*}
It's easy to see that the set $\hat{Q}$ is compact and for any $j=1,2,\dots$ the point $(g^*-\frac{1}{j},0)$ does not belong to $\hat{Q},$ where 0 is the zero element of $l_1$ (otherwise, $g^*$ is not the minimum in \eqref{M22}).  Due to Hahn-Banach separation theorem (see, e.g., \cite{DS}, Section V.2) there exists a sequence $(\k^j,\lambda^j)\in\reals\times l^\infty$ (where $\lambda^j=(\lambda_1^j,\lambda_2^j,\dots)$) such that
		\begin{equation}\label{E1}
		\begin{aligned}
			&\k^j\left(g^*-\frac{1}{j}\right)+\delta^j \leq \inf_{(\theta,x)\in\hat{Q}}\{ \k^j\theta+\sum_{i=1}^{\infty}\lambda_i^jx_i \} \\
			&=\inf_{\gamma\in\mathcal{P}(G)}\{\k^j\theta+\int_{G}(\psi_{\lambda^j}(f(y,u))-\psi_{\lambda^j}(y))\gamma(dy,du),\;
				{\rm s.t.} \; \theta\geq\int_{G}g(y,u)\gamma(dy,du)\},
			\end{aligned}
			\end{equation}
		where $\delta^j>0$ for all $j$ and $\psi_{\lambda^j}:=\sum_{i=1}^{\infty}\lambda_i^j\phi_i$. From the last formula it is easy to see that $\k^j\ge 0$. Let us show that, in fact, $\k^j>0$. Indeed, if it was not the case and $\k^j=0$, then we would have
		\begin{equation*}
		\begin{aligned}
0< \delta^j\le \min_{\gamma\in\mathcal{P}(G)}\int_{G}(\psi_{\lambda^j}(f(y,u))-\psi_{\lambda^j}(y))\gamma(dy,du)
		=	\min_{(y,u)\in G}\{(\psi_{\lambda^j}(f(y,u))-\psi_{\lambda^j}(y)),
			\end{aligned}
			\end{equation*}
	which is a contradiction to \eqref{H3}. Thus, $\k^j>0$. Dividing \eqref{E1} through by $\k_j$ we obtain
		\begin{equation*}
		\begin{aligned}
		&g^*-\frac{1}{j}<\min_{\gamma\in \P(G)}\{\int_{G}\big(g(y,u)+{1\o \k^j}(\psi_{\lambda^j}(f(y,u))-\psi_{\lambda^j}(y))\big)\gm\}\\
		&=\min_{(y,u)\in G}\{g(y,u)+{1\o \k^j}(\psi_{\lambda^j}(f(y,u))-\psi_{\lambda^j}(y))\}\le \mu^*_C(y_0).
			\end{aligned}
			\end{equation*}
Therefore, $g^*\le\mu^*_C$. Taking into account inequalities \eqref{M26} and \eqref{M24} we conclude that $g^*=\mu^*$.
\hf

\begin{Proposition}\label{P11}
The limit $\disp \lim_{\a\up 1}\min_{y\in Y}(1-\a)V_{\a}(y)$ exists and is equal to $g^*$.
\end{Proposition}
{\bf Proof.}
Let us show that
\begin{equation}\label{I1}
\limsup_{\a\up 1} \left\{\bigcup_{y_0\in Y} W_{\a}(y_0)\right\}\subset W.
\end{equation}
Indeed, let $\a_i\up 1$, $y_i\in Y$ and $\g_i\in W_{\a_i}(y_i)$ be such that $\g_i\to \g$. We have
\begin{equation*}
\begin{aligned}
0=&\int_{G}[\a_i(\ph(f(y,u))-\ph(y))+(1-\a_i)(\ph(y_i)-\ph(y))]\gamma_i(dy,du)\\
=&\int_{G}[(\a_i-1)(\ph(f(y,u))-\ph(y))+(1-\a_i)(\ph(y_i)-\ph(y))]\gamma_i(dy,du)\\
&+\int_{G}(\ph(f(y,u))-\ph(y))\gamma_i(dy,du).
\end{aligned}
\end{equation*}
Passing to the limit as $i\to \infty$ in this equality we obtain $\disp \int_{G}(\ph(f(y,u))-\ph(y))\gamma(dy,du)=0$,
therefore, $\g\in W$, i.e, \eqref{I1} holds. It follows from \eqref{I1} and \eqref{D3} that
\begin{equation}\label{I2}
g^*\le \liminf_{\a\up 1}\inf_{y\in Y} g^*_{\a}(y)=\liminf_{\a\up 1}\min_{y\in Y}(1-\a)V_{\a}(y).
\end{equation}
From \eqref{B1} it follows that for any $\a\in (0,1)$ we have
$$
0\le g(y,u)+\a V_{\a}(f(y,u))-V_{\a}(y) \quad\hbox{for all }(y,u)\in G.
$$
Therefore,
\begin{equation}\label{I7}
g(y,u)+V_{\a}(f(y,u))-V_{\a}(y)\ge (1-\a)V_{\a}(f(y,u))\ge \min_{y'\in Y}(1-\a)V_{\a}(y') \quad\hbox{for all }(y,u)\in G.
\end{equation}
Consequently, 
$$
\inf_{(y,u)\in G} \{g(y,u)+V_{\a}(f(y,u))-V_{\a}(y)\}\ge \min_{y\in Y}(1-\a)V_{\a}(y) 
$$
and
$$
\mu^*\ge \min_{y\in Y}(1-\a)V_{\a}(y).
$$ 
Along with Proposition \ref{P10}, the latter implies
$$
g^*= \mu^*\ge \limsup_{\a\up 1}\min_{y\in Y}(1-\a)V_{\a}(y).
$$
The assertion of the proposition follows from this relation and \eqref{I2}. \hf

\bigskip

The following two lemmas, proved in the Appendix, are discrete-time analogs of \cite{GruneSIAM98}, Lemma 3.5 (ii) and \cite{GruneJDE98}, Lemma 3.8.
For $v\in \reals$ the notation $[v]$ stands for the integer part of $v$.

\begin{Lemma}\label{Cor1}
Let $g:\, \T\to \reals$ be a function such that $|g(t)|\le M$ for all $t$. Let $\a\in(0,1)$ and
\begin{equation}\label{K1}
\s:=(1-\a)\sum_{t=0}^{\infty} \a^t g(t).
\end{equation}
Then for any $\ve>0$ there exists a positive integer $\disp T\ge \left[{\ve\o (4M+4|\s|+\ve)(-\ln \a)}\right]$ satisfying
\begin{equation}\label{K2}
{1\o T}\sum_{t=0}^{T-1} g(t)< \s+\ve+{2M\o T}.
\end{equation}
\end{Lemma}

\begin{Lemma}\label{L3}
Let $g:\, \T\to \reals$ be a function such that $|g(t)|\le M$ for all $t$.
Let $t$ be an arbitrary positive integer and
$$
\s:={1\o t}\sum_{\t=0}^{t-1}q(\t).
$$
For any $\ve>0$ there exists $t^*\in\{0,\dots,t-1\}$ such that
\begin{equation}\label{I21}
{1\o S}\sum_{\t=0}^{S-1}q(t^*+\t)\le \s +\ve\quad\hbox{for all } S\in\{1,\dots,t-t^*\}.
\end{equation}
Moreover,
\begin{equation}\label{K5}
l(t):=t-t^*\to \infty \quad\hbox{as }t\to \infty.
\end{equation}
\end{Lemma}

\begin{Proposition}\label{P12}
The limit $\disp \lim_{S\to \infty}G_S$ exists and is equal to $g^*$.
\end{Proposition}
{\bf Proof.} Let us show first that
\begin{equation}\label{M10}
\limsup_{S\to \infty} \G_S\subset W.
\end{equation}
Take a sequence $S_i\to \infty$ as $i\to \infty$ and let $\g_i\in \G_{S_i}$ be such that $\g_i\to\g$.
Since $\g_i\in \G_{S_i}$, there exists an initial condition $y_{0i}$ and a control $u_i(\cdot)\in \U_{S_i}(y_{0i})$ such that for the corresponding trajectory $y_i(\cdot)$ and any $\ph\in C(Y)$ we have
\begin{equation*}
\begin{aligned}
&\int_G (\ph(f(y,u))-\ph(y))\g_i(dy,du)=
{1\o S_i}\sum_{t=0}^{S_i-1} (\ph(f(y_i(t),u_i(t)))-\ph(y_i(t)))\\
&={1\o S_i}\sum_{t=0}^{S_i-1} (\ph(y_i(t+1))-\ph(y_i(t)))
={1\o S_i}(\ph(y_i(S_i)-\ph(y_{0i})).
\end{aligned}
\end{equation*}
Therefore,
\begin{equation*}
\begin{aligned}
&\int_G (\ph(f(y,u))-\ph(y))\g(dy,du)=\lim_{i\to \infty} \int_G (\ph(f(y,u))-\ph(y))\g_i(dy,du)\\
&=\lim_{i\to \infty}{1\o S_i}(\ph(y_i(S_i)-\ph(y_{0i}))=0
\end{aligned}
\end{equation*}
due to boundedness of $Y$. Thus, $\g\in W$, i.e, inclusion \eqref{M10} holds, which implies that
\begin{equation}\label{M16}
\liminf_{S\to\infty}G_S\ge g^*.
\end{equation}
Take a sequence $\a_i\up 1$. Due to Proposition \ref{P11} there exists a sequence of initial conditions $y_{0i}$, controls $u_i(\cdot)\in \U(y_{0i})$ and the corresponding trajectories $y_i(\cdot)$ such that
\begin{equation*}
(1-\a_i)\sum_{t=0}^{\infty} \a_i^t g(y_i(t),u_i(t))=g^*+\xi_i,
\end{equation*}
where $\lim_{i\to \infty} \xi_i=0$. Applying Lemma \ref{Cor1} with $\s=g^*+\xi_i$ and $\ve=\sqrt{-\ln \a_i}$  we conclude that there exists a sequence $S_i$, such that $S_i\ge K/\sqrt{-\ln \a_i}$ ($K$ is a constant independent of $i$) and
\begin{equation}\label{I25}
{1\o S_i}\sum_{t=0}^{S_i-1} g(y_i(t),u_i(t))< g^*+\xi_i+\sqrt{-\ln \a_i}+{2M\o S_i},
\end{equation}
therefore, $\disp\liminf_{S\to \infty} G_{S}\le g^*$. Together with \eqref{M16} this implies that
\begin{equation}\label{I23}
\liminf_{S\to \infty} G_{S}= g^*.
\end{equation}
The latter means that
$$
{1\o S_i}\sum_{t=0}^{S_i-1} g(y_i(t),u_i(t))=g^*+\eta_i,
$$
where $\lim_{i\to \infty}\eta_i=0$. Let us apply Lemma \ref{L3} in which $S_i$ plays the role of $t$ and $\s=g^*+\eta_i$. Set $\ve={1/{S_i}}$, denote the value corresponding to $t^*$ by $t_i$ and $l(S_i):=S_i-t_i$. We conclude that $l(S_i)\to \infty$ as $i\to \infty$ and
\begin{equation}\label{I22}
{1\o S}\sum_{t=0}^{S-1} g(y_i(t_i+t),u_i(t_i+t))\le g^*+\eta_i+{1\o {S_i}}\quad\hbox{for all }S\in \{1,\dots,l(S_i)\}.
\end{equation}
Let $\tilde u_i(\cdot)=u_i(t_i+\cdot)$, $\tilde y_i(\cdot)=y_i(t_i+\cdot)$. Note that $(\tilde u_i,\tilde y_i)$ is an admissible process. It follows from \eqref{I22} that
$$
G_S\le {1\o S}\sum_{t=0}^{S-1} g(\tilde y_i(t),\tilde u_i(t))\le g^*+\eta_i+{1\o {S_i}}\quad\hbox{for all }S\in \{1,\dots,l(S_i)\},
$$
hence,
$$
\limsup_{S\to \infty} G_S\le g^*,
$$
which, along with \eqref{I23}, completes the proof of the proposition. \hf

\bigskip

Combining the assertions of Propositions \ref{P10}, \ref{P11}, and \ref{P12}, we complete the proof of Theorem \ref{T1}.

\bigskip

The theorem below asserts convergence of the sets of occupational measures $\G_{\a}$ and $\G_{S}$ defined in Section 2 to $W$ given by \eqref{M17}.
\begin{Theorem}\label{T2}
The following holds:
$$
\lim_{\a\up 1}\rho_H (\bar{\rm co}\,\,\G_{\a},W)=\lim_{S\to \infty}\rho_H (\bar{\rm co}\,\,\G_S,W)=0.
$$
\end{Theorem}

{\bf Proof.}
The assertion of Proposition \ref{P11} in terms of occupational measures can be written as
$$
\lim_{\a\up 1}\min_{\g\in \G_{\a}}\int_G g(y,u)\g(dy,du)=\min_{\g\in W}\int_G g(y,u)\g(dy,du),
$$
which, due to linearity of the integral with respect to $\g$, implies that
$$
\lim_{\a\up 1}\min_{\g\in \bar{\rm co}\,\,\G_{\a}}\int_G g(y,u)\g(dy,du)=\min_{\g\in W}\int_G g(y,u)\g(dy,du).
$$
Since $g$ in the equality above can be any continuous function, we can write
\begin{equation}\label{M17}
\lim_{\a\up 1}\min_{\g\in \bar{\rm co}\,\,\G_{\a}}\int_G q(y,u)\g(dy,du)=\min_{\g\in W}\int_G q(y,u)\g(dy,du) \quad\hbox{for all }q\in C(G).
\end{equation}
Denote
$$
W_{\a}:=\bigcup_{y_0\in Y} W_{\a}(y_0).
$$
Due to \eqref{I1} we have
\begin{equation}\label{M20}
\limsup_{\a\up 1}W_{\a}\subset W,
\end{equation}
which, due to convexity of $W$, implies that
\begin{equation*}
\limsup_{\a\up 1}(\bar{\rm co}\,\,W_{\a})\subset W,
\end{equation*}
that is,
\begin{equation}\label{M6}
\lim_{\a\up 1}\max_{\g\in \bar{\rm co}\,\,W_{\a}}\rho(\g,W)=0.
\end{equation}
From the inclusion
$$
\G_{\a}(y_0)\subset W_{\a}(y_0) \quad\hbox{for all }y_0\in Y,
$$
proved in Proposition \ref{P4}, by taking the union with respect to $y_0\in Y$ and, then, closure of the convex hull, we conclude that
$$
\bar{\rm co}\,\,\G_{\a}\subset \bar{\rm co}\,\,W_{\a}.
$$
Therefore, from \eqref{M6} we get
\begin{equation*}
\lim_{\a\up 1}\max_{\g\in \bar{\rm co}\,\,\G_{\a}}\rho(\g,W)=0.
\end{equation*}
To complete the proof of the equality
\begin{equation}\label{M13}
\lim_{\a\up 1}\rho_H (\bar{\rm co}\,\,\G_{\a},W)=0
\end{equation}
it remains to show that
\begin{equation*}
\lim_{\a\up 1}\max_{\g\in W}\rho(\g,\bar{\rm co}\,\,\G_{\a})=0.
\end{equation*}
The proof of this relation is based on formula \eqref{M17} and weak$^*$ separation theorem. It follows the same steps
as the proof of Proposition 6.1 in \cite{GQ}, starting with formula (6.6). The only difference is that the parameter $C$, approaching 0 in \cite{GQ}, should be replaced with $\a$, approaching 1.  We do not reproduce this proof here.

The proof of the second equality of the theorem $\lim_{S\to \infty}\rho_H (\bar{\rm co}\,\,\G_S,W)=0$ is very similar to the proof of \eqref{M13}.
Namely, Proposition \ref{P12}
can be written in terms of occupational measures as
$$
\lim_{\a\up 1}\min_{\g\in \G_{S}}\int_G g(y,u)\g(dy,du)=\min_{\g\in W}\int_G g(y,u)\g(dy,du),
$$
which implies that
\begin{equation}\label{M18}
\lim_{\a\up 1}\min_{\g\in \bar{\rm co}\,\,\G_{S}}\int_G q(y,u)\g(dy,du)=\min_{\g\in W}\int_G q(y,u)\g(dy,du) \quad\hbox{for all }q\in C(G).
\end{equation}
Further, from \eqref{M10} we derive that (cf. \eqref{M20}-\eqref{M6})
\begin{equation}\label{M19}
\lim_{S\to \infty}\max_{\g\in \bar{\rm co}\,\,\G_{S}}\rho(\g,W)=0.
\end{equation}
The rest of the proof follows from \eqref{M18} and \eqref{M19} using weak$^*$ separation theorem following the lines of \cite{GQ}, as described above.

\hf

\section{Appendix}\label{Sec-Appendix}

{\bf Proof of Proposition \ref{P3}.} We have
\begin{equation*}
\begin{aligned}
V(y_0)&=\min_{u(\cdot)\in \U(y_0)} \sum_{t=0}^{\infty} \a^t g(y(t),u(t))
=\min_{u(\cdot)\in \U(y_0)}\{g(y(0),u(0))+\sum_{t=1}^{\infty} \a^t g(y(t),u(t))\}\\
&=\min_{u(0)\in A(y(0))}\{g(y(0),u(0))+\a\min_{\{u(t)\in A(y(t)),\,t\ge 1\}} \sum_{t=1}^{\infty} \a^{t-1} g(y(t),u(t))\}.
\end{aligned}
\end{equation*}
The second minimum is equal to $V(y(1))=V(f(y(0),u(0)))$, therefore,
$$
V(y_0)=\min_{u(0)\in A(y(0))}\{g(y(0),u(0))+\a V(f(y(0),u(0)))\}.
$$
Replacing now $u(0)$ and $y(0)$ with $u$ and $y$, respectively, we obtain relation \eqref{B1}.
\hfill{$\Box$}

\bigskip

\begin{Lemma}\label{Gr1}{\rm (\cite{GruneSIAM98}, Lemma 3.5 (ii))}
Let $q:\, [0,\infty)\to \reals$ be a measurable function such that $|q(\t)|\le M$ for a.a. $\t\in \reals$. Let $\d>0$ be arbitrary and
\begin{equation}\label{K11}
\ts:=\d\int_0^{\infty} e^{-\d \t}q(\t)\,d\t.
\end{equation}
Then for any $\ve>0$ there exists  $\disp \tilde T\ge {\ve\o (4M+4|\ts|+\ve)\d}$ satisfying
\begin{equation}\label{K3}
{1\o \tilde T}\int_0^{\tilde T} q(\t)\,d\t\le \ts+\ve.
\end{equation}
\end{Lemma}
{\bf Proof of Lemma \ref{Cor1}.} Lemma \ref{Cor1} is a discrete-time analog of Lemma \ref{Gr1}.

Define the piecewise constant function $q:\, [0,\infty)\to \reals$ by
\begin{equation*}
q(\t)=g(t), \quad \t\in[t,t+1),\,t\in\T
\end{equation*}
and apply Lemma \ref{Gr1} with $\d=-\ln\a$.  Let us first evaluate $\ts$ given by \eqref{K11}. For $t\in \T$ we have
$$
\int_{t}^{t+1} e^{-\d \t}q(\t)\,d\t={1\o \d}(1-e^{-\d})g(t)e^{-\d t},
$$
therefore,
$$
\ts:=\d\int_{0}^{\infty} e^{-\d \t}q(\t)\,d\t=(1-e^{-\d})\sum_{t=0}^{\infty}e^{-\d t}g(t)=(1-\a)\sum_{t=0}^{\infty} \a^t g(t)=\s.
$$
Due to Lemma \ref{Gr1} there exists $\tilde T\ge {\ve/\big( (4M+4|\s|+\ve)(-\ln \a)\big)}$ such that
\begin{equation}\label{M2}
{1\o \tilde T}\int_0^{\tilde T} q(\t)\,d\t\le \s+\ve.
\end{equation}
In the case if $0<\tilde T<1$, then ${1\o \tilde T}\int_0^{\tilde T} q(\t)\,d\t=g(0)$ and inequality \eqref{K2} holds in the form
\begin{equation*}
{1\o T}\sum_{t=0}^{T-1} g(t)\le \s+\ve
\end{equation*}
with $T=1$. Assume, therefore, that $\tilde T\ge 1$.

Let $T:=[\tilde T]\ge 1$ and denote $\D T:=\tilde T-T\in [0,1)$. We have
\begin{equation}\label{M1}
{1\o \tilde T}\int_0^{\tilde T} q(\t)\,d\t= {1\o T+\D T}\int_0^{T} q(\t)\,d\t+{1\o T+\D T}\int_{T}^{T+\D T} q(\t)\,d\t.
\end{equation}
For the second integral we have
\begin{equation}\label{M3}
{1\o T+\D T}\int_{T}^{T+\D T} q(\t)\,d\t\ge -{M\D T\o T+\D T} > -{M\o T}.
\end{equation}
Taking into account that $1/(1+x)\ge 1-x$ for $x> -1$ we have
$$
{1\o T+\D T}={1\o T}{1\o 1+\D T/ T}\ge {1\o T}-{\D T\o T^2},
$$
therefore, in the case if $\int_0^T q(\t)\,d\t\ge 0$, for the first integral on the right hand side of \eqref{M1} we have
\begin{equation}\label{M4}
{1\o T+\D T}\int_0^{T} q(\t)\,d\t \ge{1\o T}\int_0^{T} q(\t)\,d\t-{M \D T\o T}> {1\o T}\int_0^{T} q(\t)\,d\t-{M \o T}.
\end{equation}
If  $\int_0^T q(\t)\,d\t< 0$, then $\disp{1\o T+\D T}\int_0^{T} q(\t)\,d\t \ge {1\o T}\int_0^{T} q(\t)\,d\t$ and the inequality above still holds.
Thus, we obtain from \eqref{M1}-\eqref{M4},  that
\begin{equation}\label{M11}
{1\o \tilde T}\int_0^{\tilde T} q(\t)\,d\t> {1\o T}\int_0^{T} q(\t)\,d\t-{2M\o T}={1\o T}\sum_{t=0}^{T-1} g(t)-{2M\o T},
\end{equation}
and \eqref{K2} follows from \eqref{M2} and \eqref{M11}. \hf

\bigskip

{\bf Proof of Lemma \ref{L3}.}
Let $\disp\b:=\max_{1\le s\le t} {1\o s}\sum_{\t=0}^{s-1}q(\t).$ If $\b\le \s+\ve$ then the statement of the lemma holds with $t^*=0$. Assume, therefore, that
$\b>\s+\ve$ and set
$$
t^*=\max\{s\in\{1,\dots,t\}|\,{1\o s}\sum_{\t=0}^{s-1}q(\t)> \s+\ve\}.
$$
Let us show that this $t^*$ satisfies the required properties.
Indeed, $t^*\neq t$ due to the definition of $\s$, hence, $0\le t^*\le t-1$. Let us show that \eqref{I21} is satisfied. Assume the contrary, that is, there exists $1\le s_1\le t-t^*$ such that
$\disp \s+\ve<{1\o s_1}\sum_{\t=0}^{s_1-1}q(t^*+\t)={1\o s_1}\sum_{\t=t^*}^{t^*+s_1-1}q(\t)$. This implies that
$$
\sum_{\t=0}^{t^*+s_1-1}q(\t)=\sum_{\t=0}^{t^*-1}q(\t)+\sum_{\t=t^*}^{t^*+s_1-1}q(\t)>(\s+\ve)t^*+(\s+\ve)s_1=(\s+\ve)(t^*+s_1),
$$
which contradicts the definition of $t^*$.

Let us show now that $l(t):=t-t^*\to \infty$ as $t\to \infty$. We have
$$
\s t=\sum_{\t=0}^{t-1}q(\t)=\sum_{\t=0}^{t^*-1}q(\t)+\sum_{\t=t^*}^{t-1}q(\t)>(\s+\ve)t^*-(t-t^*)M.
$$
This can be equivalently written as
$$
\s t>(\s+\ve)(t-l)-lM,
$$
or,
$$
l(\s+\ve+M)> \ve t,
$$
which implies that $l\to \infty$ as $t\to \infty$, that is, \eqref{K5} holds.
\hf

\bigskip
\bigskip

Email addresses of the authors: 

V. Gaitsgory vladimir.gaitsgory@mq.edu.au\\
A. Parkinson alex.parkinson@students.mq.edu.au\\
I. Shvartsman ius13@psu.edu (corresponding author)

\end{document}